\documentclass{birkjour_t2}
\usepackage{amsmath}
\usepackage{tkz-euclide}
\usetkzobj{all}
\tikzset{elegant/.style={smooth,thick,samples=50,cyan}}
\tikzset{eaxis/.style={->,>=stealth}}
\usepackage{color,tikz}
\usetikzlibrary{decorations.pathreplacing,calc}
\usepackage{rotating}
\tikzset{liltext/.style={font=\tiny}}
\usepackage{amsfonts}
\usepackage{amssymb}
\usepackage[english]{babel}
\usepackage{color}
\usepackage{graphicx}
\usepackage{lmodern}
\usepackage{amsthm}
\usepackage{mathrsfs}
\usepackage{microtype}
\usepackage{mathscinet}
\usepackage{enumitem}
\usepackage[cal=boondoxo,bb=ams]{mathalfa}
\usepackage{hyperref}
\hypersetup{hidelinks}

\renewcommand{\theequation}{\arabic{section}.\arabic{equation}}

\newtheorem{defn}{Definition}[section]
\newtheorem{thm}{Theorem}[section]

\newtheorem{rem}{Remark}[section]

\newcommand{\ml}{\mathcal}
\newcommand{\mb}{\mathbb}

\begin{document}

%-------------------------------------------------------------------------
% editorial commands: to be inserted by the editorial office
%
%\firstpage{1} \volume{228} \Copyrightyear{2004} \DOI{003-0001}
%\seriesextraline{This is the Concrete Title of this Book\br H.E. R and S.T.C. W, Eds.}
%\firstpage{1}
%      environment-name

\title[Semilinear Moore -- Gibson -- Thompson equation with nonlinearity of derivative type]
 {A blow -- up result for the semilinear Moore -- Gibson -- Thompson equation with nonlinearity of derivative type in the conservative case}

%----------Author 1
\author[W. Chen]{Wenhui Chen}
\address{Institute of Applied Analysis, Faculty of Mathematics and Computer Science\\
	 Technical University Bergakademie Freiberg\\
	  Pr\"{u}ferstra{\ss}e 9\\
	   09596 Freiberg\\
	    Germany}
\email{wenhui.chen.math@gmail.com}

%----------Author 2
\author[A. Palmieri]{Alessandro Palmieri}

\address{Department of Mathematics\\
	 University of Pisa\\
	Largo B. Pontecorvo 5\\
	56127 Pisa\\
	Italy}
\email{alessandro.palmieri.math@gmail.com}

\subjclass{Primary  35B44, 35L30, 35L76 ; Secondary 35B33, 35L25}
%\MSC[2010] Primary 35L52, 35L99; Secondary 35B33, 35B44

\keywords{{Moore -- Gibson -- Thompson equation, semilinear third order equation, blow-up, lifespan estimates, Glassey exponent.}
}
\date{June 10, 2019}

\begin{abstract} 
	In this paper, we study  the blow -- up of solutions to the semilinear Moore -- Gibson -- Thompson (MGT) equation with nonlinearity of derivative type $|u_t|^p$ in the conservative case. We apply an iteration method in order to study both the subcritical case and the critical case. Hence, we obtain a blow -- up result  for the semilinear MGT equation (under suitable assumptions for initial data) when the exponent $p$ for the nonlinear term  satisfies $1<p\leqslant (n+1)/(n-1)$  for $n\geqslant2$ and $p>1$ for $n=1$. In particular, we find the same blow -- up range for $p$ as in the corresponding semilinear wave equation with nonlinearity of derivative type.
\end{abstract}

%%% ----------------------------------------------------------------------
\maketitle
%%% ----------------------------------------------------------------------
%\tableofcontents
\section{Introduction}

Over the last years, the Moore -- Gibson -- Thompson (MGT) equation  (cf. \cite{MooreGibson1960,Thompson1972}), a linearization of a model for the wave propagation in viscous thermally relaxing fluids, has been studied by several authors %\textcolor{red}{$\blacktriangleright$ old version} (see \cite{Stokes1851,MooreGibson1960,Thompson1972,LebonCloot1989,NaugolnykhOstrovsky1998,Gorain2010,KaltenbacherLasieckaMarchand2011,KaltenbacherLasiecka2012,MarchandMcDevittTriggiani2012,Jordan2014,LasieckaWang2015,PellicerSola2015,CaixetaLasieckaDominos2016,DellOroLasieckaPata2016,LasieckaWang2016,DellOroPata2017,Lasiecka2017,PellicerSaiHouari2017,AlvesCaixetaSilvaRodrigues2018,DellOroLasieckaPata2019}).  \textcolor{red}{new version} (see, for example, \cite{Gorain2010,KaltenbacherLasieckaMarchand2011,KaltenbacherLasiecka2012,MarchandMcDevittTriggiani2012,Jordan2014,LasieckaWang2015,CaixetaLasieckaDominos2016,DellOroLasieckaPata2016,LasieckaWang2016,DellOroPata2017,Lasiecka2017,PellicerSaiHouari2017,AlvesCaixetaSilvaRodrigues2018,DellOroLasieckaPata2019,PellicerSola2015}).  
%\textcolor{red}{$\blacktriangleleft$} 
(see, for example, \cite{Gorain2010,KaltenbacherLasieckaMarchand2011,KaltenbacherLasiecka2012,MarchandMcDevittTriggiani2012,Jordan2014,LasieckaWang2015,CaixetaLasieckaDominos2016,DellOroLasieckaPata2016,LasieckaWang2016,DellOroPata2017,Lasiecka2017,PellicerSaiHouari2017,AlvesCaixetaSilvaRodrigues2018,DellOroLasieckaPata2019,PellicerSola2015}).

%\textcolor{red}{$\blacksquare$ See in the bibliography the comments on \cite{Stokes1851,LebonCloot1989,NaugolnykhOstrovsky1998}} \qquad\textcolor{purple}{AGREE !}

%\textcolor{red}{$\blacksquare$ Please check that  \cite{MooreGibson1960,Thompson1972} are the works due to which the equation we are considering is named MGT, I could not find them!} \qquad\textcolor{purple}{AGREE ! I CANNOT FIND WHO GIVE THE NAME OF THE EQUATION MGT !}

 This model is realized through the third order hyperbolic equation
\begin{equation}\label{General.MGT.Equation}
\tau u_{ttt}+ u_{tt}-c^2\Delta u-b\Delta u_t=0.
\end{equation} 
In the physical context of acoustic waves, the unknown function $u=u(t,x)$ denotes the scalar acoustic velocity, $c$ denotes the speed of sound and $\tau$ denotes the thermal relaxation. Besides, the coefficient $b=\beta c^2$ is related to the diffusivity of the sound with $\beta\in (0,\tau]$.  Let us point out that there is the transition from a linear model  that can be described in the case of bounded domains with an exponentially stable strongly continuous semigroup in the case $0<\beta<\tau$ to the limit case $\beta=\tau$, where the exponential stability of a semigroup is lost and it holds the conservation of a suitable defined energy (see \cite{KaltenbacherLasieckaMarchand2011,MarchandMcDevittTriggiani2012}). For this reason, we will refer to the limit case $\beta=\tau$ as to the \emph{conservative case} throughout this paper.

We consider the semilinear Cauchy problem for MGT equation in the conservative case with nonlinearity of derivative type, namely,
\begin{equation}\label{Semi.MGT.u_t.epsilon}
\begin{cases}
\beta u_{ttt}+u_{tt}-\Delta u-\beta \Delta u_t=|u_t|^p, & x\in\mb{R}^n,\, t>0, \\
(u,u_t,u_{tt})(0,x)= \varepsilon ( u_0, u_1, u_2)(x), &x\in\mb{R}^n,
\end{cases}
\end{equation}
 where $p>1$  and $\varepsilon$ is a positive parameter describing the size of initial data. Note that, for the sake of simplicity, we normalized the speed of the sound by putting $c^2=1$. We are interested in investigating the blow -- up in finite time of local (in time) solutions under suitable sign assumptions for the Cauchy data regardless of their size.  Let us underline that, while the MGT equation has been widely investigated in the case of bounded domains via semigroups theory, very few results concerning nonlinear Cauchy problems for MGT equation are available up to the knowledge of the authors. In \cite{Racke2019}, the semilinear Cauchy problem with nonlinearity $\partial_t\left(k (u_t)^2+|\nabla u|^2\right)$ is considered in the \emph{dissipative case} $0<\tau<\beta$, where $k$ is a suitable  constant. More precisely, a global existence result for small data solutions is proved providing that initial data are sufficiently regular and satisfy certain integral relations (cf. \cite[Theorem 5.1]{Racke2019}). Moreover, a blow -- up result for the conservative case with power nonlinearity can be found in \cite{ChenPal19}. 

Let us provide some results which are related to our model \eqref{Semi.MGT.u_t.epsilon}. By choosing $\beta=0$, we find that \eqref{Semi.MGT.u_t.epsilon} corresponds formally to the semilinear wave equation
\begin{equation}\label{Semi.Wave.Eq}
\begin{cases}
u_{tt}-\Delta u=|u_t|^p,& x\in\mb{R}^n,\,t>0,\\
(u,u_t)(0,x)=\varepsilon(u_0,u_1)(x),&x\in\mb{R}^n,
\end{cases}
\end{equation}
where $p>1$. According to \cite{John1981,Sideris1983,Masuda1983,Schaeffer1986,Rammaha1987,Agemi1991,HidanoTsutaya1995,Tzvetkov1998,Zhou2001,HidanoWangYokoyama2012}  the critical exponent of \eqref{Semi.Wave.Eq} is the so -- called \emph{Glassey exponent} $p_{\mathrm{Gla}}(n) \doteq (n+1)/(n-1)$. Moreover, the sharp behavior of the lifespan $T(\varepsilon)$ of local (in time) solutions to \eqref{Semi.Wave.Eq} with respect to a sufficiently small parameter $\varepsilon>0$ is given by
\begin{align*}
T(\varepsilon)\approx
	\begin{cases}
C \varepsilon^{-(\frac{1}{p-1} -\frac{n-1}{2})^{-1}}&\mbox{if}\ \ 1<p<p_{\mathrm{Gla}}(n),\\
\exp\left(C\varepsilon^{-(p-1)}\right)&\mbox{if}\ \ p=p_{\mathrm{Gla}}(n).
\end{cases}
\end{align*}

The main result of this paper consists of a blow -- up result for \eqref{Semi.MGT.u_t.epsilon} when the power of the nonlinear term is in the sub -- Glassey range (including the case $p=p_{\mathrm{Gla}}(n) $).

In order to prove this result, we are going to apply an iterative argument for a suitable time -- dependent functional, which depends on a local (in time) solution to \eqref{Semi.MGT.u_t.epsilon}. For the choice of the functional we follow \cite{LaiTakamura19} whereas concerning the iteration procedure we use some key ideas from \cite{ChenPal19}, where a technique to deal with an unbounded exponential multiplier in the iteration frame is developed. This approach is based on the idea of slicing the interval of integration and it has been introduced by Takamura and coauthors in the study of critical cases for wave models (see \cite{AKT00,TakWak11,TakWak14,WakYor18} for example). Recently, many papers have been devoted to the study of blow -- up results for semilinear second order hyperbolic models with the aid of a time dependent multipliers. The first paper in this direction is \cite{LT18Scatt} followed then by \cite{LaiTakamura19,LT18ComNon,PalTak19,PalTak19dt,PalTak19mix,ITW19}. In these papers, the time -- dependent multiplier is bounded by positive constants from above and from below and it is used to study semilinear damped wave models with time -- dependent coefficients for the damping terms in the scattering producing case. On the other hand, the case of unbounded time -- dependent multipliers is considered for semilinear wave models with scale -- invariant damping and mass terms in \cite{LTW17,LaiTakamura19,PT18,Pal19,PT19}. 

Before stating the main result of this paper, let us introduce a suitable notion of energy solutions to the Cauchy problem \eqref{Semi.MGT.u_t.epsilon}.
\begin{defn}\label{Defn.Energy.Solution}
	Let $(u_0,u_1,u_2)\in H^2(\mb{R}^n)\times H^1(\mb{R}^n)\times L^2(\mb{R}^n)$. We say that $u$ is an energy solution of \eqref{Semi.MGT.u_t.epsilon} on $[0,T)$ if
	\begin{align*}
	u\in\ml{C}([0,T),H^2(\mb{R}^n))\cap \ml{C}^1([0,T),H^1(\mb{R}^n))\cap \ml{C}^2([0,T),L^2(\mb{R}^n))\quad\text{such that}\quad u_t\in L_{\emph{loc}}^p([0,T)\times\mb{R}^n)
	\end{align*}
	satisfies $u(0,\cdot)=\varepsilon u_0$ in $H^2(\mb{R}^n)$ and the integral relation
	\begin{align}\label{Def.Energy solution u_t}
	\beta\int_{\mb{R}^n}&u_{tt}(t,x)\phi(t,x) \,\mathrm{d}x+\int_{\mb{R}^n}u_t(t,x)\phi(t,x)\,\mathrm{d}x-\beta\varepsilon\int_{\mb{R}^n} u_2(x)\phi(0,x)\,\mathrm{d}x-\varepsilon\int_{\mb{R}^n} u_1(x)\phi(0,x)\,\mathrm{d}x\notag\\
	&+\beta\int_0^t\int_{\mb{R}^n}\left(\nabla u_t(s,x)\cdot\nabla\phi(s,x)-u_{tt}(s,x)\phi_t(s,x)\right)\mathrm{d}x\,\mathrm{d}s\notag\\
	&+\int_0^t\int_{\mb{R}^n}\left(\nabla u(s,x)\cdot\nabla\phi(s,x)-u_t(s,x)\phi_t(s,x)\right)\mathrm{d}x\,\mathrm{d}s\notag\\
	&=\int_0^t\int_{\mb{R}^n}|u_t(s,x)|^p\phi(s,x)\,\mathrm{d}x\,\mathrm{d}s
	\end{align}
	for any $\phi\in\ml{C}_0^{\infty}([0,T)\times\mb{R}^n)$ and any $t\in(0,T)$.
\end{defn}

Applying further steps of integration by parts in \eqref{Def.Energy solution u_t}, we get
\begin{align*}
\beta\int_{\mb{R}^n}&\big(u_{tt}(t,x)\phi(t,x)-u_t(t,x)\phi_t(t,x)-u(t,x)\Delta\phi(t,x)+u(t,x)\phi_{tt}(t,x)\big)\mathrm{d}x\\
&+\int_{\mb{R}^n}\big(u_t(t,x)\phi(t,x)-u(t,x)\phi_t(t,x)\big)\mathrm{d}x\\
&-\beta \varepsilon\int_{\mb{R}^n}\big(u_2(x)\phi(0,x)- u_1(x)\phi_t(0,x)-u_0(x)\Delta\phi(0,x)+ u_0(x)\phi_{tt}(0,x)\big)\mathrm{d}x\\
&- \varepsilon\int_{\mb{R}^n}\big(u_1(x)\phi(0,x)- u_0(x)\phi_t(0,x)\big)\mathrm{d}x\\
&+\int_0^t\int_{\mb{R}^n}u(s,x)\big(-\beta\phi_{ttt}(s,x)+\phi_{tt}(s,x)-\Delta\phi(s,x)+\beta\Delta\phi_t(s,x)\big)\mathrm{d}x\,\mathrm{d}s\\
&=\int_0^t\int_{\mb{R}^n}|u_t(s,x)|^p\phi(s,x)\, \mathrm{d}x\,\mathrm{d}s.
\end{align*}
Letting $t\rightarrow T$, we find that $u$ fulfills the definition of weak solution to \eqref{Semi.MGT.u_t.epsilon}.

%\noindent\textcolor{blue}{$\blacktriangle$ Do you think we should show something about local existence with support condition? However, I think that it is a little boring to follow the proof of our PART I. It seems that we only need to modify few steps.}
%\textcolor{red}{$\blacktriangle$ No, I think it is not necessary since it is boring and space-consuming. We could just cite our first paper and the corresponding result and write something like ``modifying slightly the proof of Theorem/Proposition ... one can prove the existence of local in time energy solutions with support contained in the forward cone ... for any $1<p\leqslant ...$ ''} \qquad\textcolor{purple}{AGREE !}

We now state our main result.
\begin{thm} \label{Thm blow-up}
	Let us consider $p>1$ such that 
	\begin{align*}
	\begin{cases} p<  \infty & \mbox{if} \ \  n=1, \\ p\leqslant  p_{\mathrm{Gla}}(n) & \mbox{if} \ \ n\geqslant 2. \end{cases}
	\end{align*}
	Let $(u_0,u_1,u_2)\in H^2(\mathbb{R}^n)\times H^1(\mathbb{R}^n)\times L^2(\mathbb{R}^n)$ be nonnegative and compactly supported functions with supports contained in $B_R$ for some $R>0$ such that $u_1$ or $u_2$ is not identically zero.\\
	Let $u$ be the energy solution to the Cauchy problem \eqref{Semi.MGT.u_t.epsilon} with lifespan $T(\varepsilon)$ satisfying
	\begin{align*}
	\mathrm{supp}\, u(t,\cdot) \subset B_{t+R} \quad \mbox{for any} \ t\in (0,T).
	\end{align*} Then, there exists a positive constant $\varepsilon_0=\varepsilon_0(u_0,u_1,u_2,n,p,R,\beta)$ such that for any $\varepsilon\in(0,\varepsilon_0]$ the solution $u$ blows up in finite time. Furthermore, the upper bound estimate for the lifespan
	\begin{align*}
	T(\varepsilon)\leqslant
	\begin{cases}
	C \varepsilon^{-(\frac{1}{p-1} -\frac{n-1}{2})^{-1}}&\mbox{if}\ \ 1<p<p_{\mathrm{Gla}}(n),\\
	\exp\left(C\varepsilon^{-(p-1)}\right)&\mbox{if}\ \ p=p_{\mathrm{Gla}}(n),\\
	\end{cases}
	\end{align*} holds, where $C>0$ is an independent of $\varepsilon$ constant.
\end{thm}

\begin{rem}  We point out that the solution  to the linear Cauchy problem for MGT equation 
\begin{equation}\label{LinearMGT.epsilon}
\begin{cases}
\beta u_{ttt}+u_{tt}-\Delta u-\beta \Delta u_t=F(t,x), & x\in\mb{R}^n,\, t>0, \\
(u,u_t,u_{tt})(0,x)= \varepsilon ( u_0, u_1, u_2)(x), &x\in\mb{R}^n,
\end{cases}
\end{equation}
 fulfills the inhomogeneous wave equation
\begin{equation}\label{MGT.Equation.Wave}
\begin{cases}
u_{tt}-\Delta u=\varepsilon\,\mathrm{e}^{-t/\beta}(u_2(x)-\Delta u_0(x))+1/\beta \int_0^t \mathrm{e}^{(\tau-t)/\beta} F(\tau,x) \, \mathrm{d}\tau, &x\in\mb{R}^n,\,t>0,\\
(u,u_t)(0,x)=\varepsilon(u_0,u_1)(x),&x\in\mb{R}^n.
\end{cases}
\end{equation}
Thus, we claim that $\mathrm{supp } \,u(t,\cdot)\subset B_{t+R}$, if we assume for some $R>0$ that $\mathrm{supp } \, u_j\subset B_R$ for any $j=0,1,2$  and $\mathrm{supp } \,F(t,\cdot)\subset B_{t+R}$ for any $t\geqslant 0$. Indeed, the source term $$f(t,x)=\varepsilon \,\mathrm{e}^{-t/\beta}(u_2(x)-\Delta u_0(x))+1/\beta \int_0^t \mathrm{e}^{(\tau-t)/\beta} F(\tau,x) \, \mathrm{d}\tau$$ in \eqref{MGT.Equation.Wave} has support contained in the forward cone $\{(t,x): |x|\leqslant t+R\}$ under these assumptions and we can use the property of finite speed of propagation for the classical wave equation. Therefore, the support condition in Theorem \ref{Thm blow-up} for a local in time solution to \eqref{Semi.MGT.u_t.epsilon} is meaningful.
\end{rem}

\noindent\textbf{Notation: } We give some notations to be used in this paper. We write $f\lesssim g$ when there exists a positive constant $C$ such that $f\leqslant Cg$. Moreover, we write $g\lesssim f \lesssim g$ by $f\approx g$. $B_R$ denotes the ball around the origin with radius $R$ in $\mathbb{R}^n$. Finally, as in the introduction, $p_{\mathrm{Gla}}(n)$ denotes the Glassey exponent.

\section{Blow -- up result in the subcritical case}\label{Section.Blow-up.Subcritical}
\subsection{Iteration frame} 
Let us consider the eigenfunction $\Phi$ of the Laplace operator  on the whole space 
\begin{align*}
\Phi(x) & \doteq
  \mathrm{e}^{x}+\mathrm{e}^{-x}  \qquad  \qquad \ \mbox{if} \ n=1, \\  \Phi(x) & \doteq \int_{\mathbb{S}^{n-1}} \mathrm{e}^{x\cdot \omega} \, \mathrm{d} \sigma_\omega \qquad \mbox{if} \ n\geqslant 2,
\end{align*} for any $x\in \mathbb{R}^n$. This function has been employed in the study of blow -- up results for the semilinear wave model in the critical case in \cite{YordanovZhang2006}. The function $\Phi$ is positive and smooth and satisfies the following remarkable properties:
\begin{align}
& \Delta \Phi =\Phi , \label{Laplace Phi =Phi} \\
& \Phi (x) \sim |x|^{-\frac{n-1}{2}}  \mathrm{e}^x \qquad \mbox{as} \ |x|\to \infty. \label{Asymptotic behavior Phi}
\end{align} 

Hence, we define the function with separate variables $\Psi=\Psi(t,x)\doteq \mathrm{e}^{-t}\Phi(x)$. Therefore, $\Psi$ is a solution of the adjoint equation to the homogeneous linear MGT equation, namely, 
\begin{align} \label{adjoint MGT eq}
-\beta\, \partial_t^3 \Psi+\partial_t^2 \Psi -\Delta \Psi +\beta \Delta\partial_t \Psi=0.
\end{align} 
By using the asymptotic behavior of $\Psi$ (cf. \cite[Equation (3.5)]{LaiTakamura19}), it follows that there exists a constant $C_1=C_1(n,R)>0$ such that
\begin{align} \label{L1 norm Psi B t+R}
\int_{B_{t+R}}\Psi(t,x)\mathrm{d}x\leqslant C_1(t+R)^{(n-1)/2}\quad\text{for any }t\geqslant0.
\end{align}

%\color{purple}
Moreover, modifying slightly the proof of Theorem 3.1 in \cite{ChenPal19} one can prove the existence of local in time energy solutions with support contained in the forward cone $\{(t,x)\in [0,T)\times\mathbb{R}^n: |x|\leqslant t+R\}$ for any $p>1$ such that $p\leqslant n/(n-2)$ when $n\geqslant3$, if $(u_0,u_1,u_2)\in H^2(\mb{R}^n)\times H^1(\mb{R}^n)\times L^2(\mb{R}^n)$ are compactly supported functions with supports contained in $B_R$ for some $R>0$.
%\color{black}

Since $u$ is supported in a forward cone, we may apply the definition of energy solution even though the test function is not compactly supported. So, applying the definition of energy solution with $\Psi$ as test function in \eqref{Def.Energy solution u_t}, we get for any $t\in(0,T)$
\begin{align}\label{Psi.Calculationu_t}
\int_0^t\int_{\mb{R}^n}|u_t(s,x)|^p\Psi(s,x)\mathrm{d}x\,\mathrm{d}s&=\beta\int_{\mb{R}^n}u_{tt}(t,x)\Psi(t,x)\mathrm{d}x+\int_{\mb{R}^n}u_t(t,x)\Psi(t,x)\mathrm{d}x\notag\\
&\quad-\beta \varepsilon\int_{\mb{R}^n} u_2(x)\Phi(x)\mathrm{d}x- \varepsilon\int_{\mb{R}^n} u_1(x)\Phi(x)\mathrm{d}x\notag\\
&\quad+\beta\int_0^t\int_{\mb{R}^n}\big(\nabla u_t(s,x)\cdot\nabla\Psi(s,x)-u_{tt}(s,x)\Psi_t(s,x)\big)\mathrm{d}x\,\mathrm{d}s\notag\\
&\quad+\int_0^t\int_{\mb{R}^n}\big(\nabla u(s,x)\cdot\nabla\Psi(s,x)-u_t(s,x)\Psi_t(s,x)\big)\mathrm{d}x\,\mathrm{d}s.
\end{align}
Consequently, performing integration by parts in \eqref{Psi.Calculationu_t} and employing the properties of $\Psi$, we find
\begin{align}\label{Eq.01}
\int_0^t \int_{\mathbb{R}^n} |u_t(s,x)|^p \Psi(s,x) \mathrm{d}x \, \mathrm{d}s&=\int_{\mathbb{R}^n} \big(\beta u_{tt}(t,x)+ (\beta+1)u_t(t,x)+u(t,x)\big)\Psi(t,x) \,\mathrm{d}x\notag\\
&\quad-\varepsilon\int_{\mathbb{R}^n} \big(\beta u_{2}(x)+ (\beta+1)u_1(x)+u_0(x)\big)\Phi(x) \,\mathrm{d}x.
\end{align}

Let us introduce
\begin{align*}
I_{\beta}[u_0,u_1,u_2]&\doteq\int_{\mathbb{R}^n} \big(\beta u_{2}(x)+ (\beta+1)u_1(x)+u_0(x)\big)\Phi(x) \,\mathrm{d}x,\\
F_1(t)&\doteq\int_{\mb{R}^n}u_t(t,x)\Psi(t,x)\,\mathrm{d}x.
\end{align*}  The functional $F_1$ will play a central role in the iteration argument, as it is the time -- dependent quantity that blows up in finite time or, in other words, it is the function that will be estimated from below iteratively.
By using these notations, we may rewrite \eqref{Eq.01} as
\begin{align}\label{Eq.02}
&\beta F'_1(t)+(2\beta+1)F_1(t)+\int_{\mb{R}^n}u(t,x)\Psi(t,x)\,\mathrm{d}x% \notag\\ &
=\int_0^t \int_{\mathbb{R}^n} |u_t(s,x)|^p \Psi(s,x)\, \mathrm{d}x \, \mathrm{d}s +\varepsilon I_{\beta}[u_0,u_1,u_2].
\end{align}
Furthermore,  the differentiation of \eqref{Eq.02} with respect to $t$ provides
\begin{align}\label{Eq.03}
\beta F''_1(t)+(2\beta+1)F'_1(t)+F_1(t)-\int_{\mb{R}^n}u(t,x)\Psi(t,x)\,\mathrm{d}x=& \int_{\mathbb{R}^n} |u_t(t,x)|^p \Psi(t,x)\, \mathrm{d}x.
\end{align}
Adding up \eqref{Eq.02} with \eqref{Eq.03}, we immediately obtain
\begin{align}\label{Eq.04}
&\beta F''_1(t)+(3\beta+1)F'_1(t)+(2\beta+2)F_1(t)\notag\\
&=\int_0^t \int_{\mathbb{R}^n} |u_t(s,x)|^p \Psi(s,x)\, \mathrm{d}x \, \mathrm{d}s+\int_{\mathbb{R}^n} |u_t(t,x)|^p \Psi(t,x) \,\mathrm{d}x+\varepsilon I_{\beta}[u_0,u_1,u_2].
\end{align}

Next, let us set
\begin{align*}
G(t)&\doteq F_1'(t)+2F_1(t)-(\beta+1)^{-1}\int_0^t\int_{\mb{R}^n}|u_t(s,x)|^p\Psi(s,x)\,\mathrm{d}x\,\mathrm{d}s-\varepsilon (\beta+1)^{-1} J_{\beta}[u_1,u_2],
\end{align*} where
\begin{align*}
J_{\beta}[u_1,u_2]&\doteq\int_{\mathbb{R}^n} \big(\beta u_{2}(x)+ (\beta+1)u_1(x)\big)\Phi(x) \, \mathrm{d}x.
\end{align*}   The auxiliary functional $G$, together with $H$ whose definition is going to be introduced in few lines, is important to derive a first lower bound estimate for $F_1$ and the iteration frame for $F_1$.
Employing \eqref{Eq.04} and the nonnegativity of $u_0$, we arrive at
\begin{align*}
\beta G'(t)+(\beta+1)G(t)=(\beta+1)^{-1} \int_{\mb{R}^n}|u_t(t,x)|^p\Psi(t,x)\,\mathrm{d}x+\varepsilon\int_{\mb{R}^n}u_0(x)\Phi(x)\,\mathrm{d}x\geqslant0,
\end{align*}
which implies in turn
\begin{align*}
G(t)\geqslant \mathrm{e}^{-(1+1/\beta)t}G(0)=\varepsilon (\beta+1)^{-1} \mathrm{e}^{-(1+1/\beta)t}\int_{\mb{R}^n}u_2(x)\Phi(x)\,\mathrm{d}x\geqslant 0,
\end{align*}
where we used the nonnegativity of $u_2$.

Combining the definition of $G$ into the inequality $G(t)\geqslant 0$, we get
\begin{align}\label{Eq.05}
F_1'(t)+2F_1(t)\geqslant (\beta+1)^{-1}\int_0^t\int_{\mb{R}^n}|u_t(s,x)|^p\Psi(s,x)\,\mathrm{d}x\,\mathrm{d}s+\varepsilon (\beta+1)^{-1} J_{\beta}[u_1,u_2] \doteq H(t) .
\end{align}
This leads to
\begin{align}\label{Initial.Lower.Bound}
F_1(t)&\geqslant \mathrm{e}^{-2t}F_1(0)+\frac{\varepsilon}{2(\beta+1)}J_{\beta}[u_1,u_2]\left(1-\mathrm{e}^{-2t}\right)\notag\\
& =\frac{\varepsilon}{2}\left(1+\mathrm{e}^{-2t}\right)\int_{\mb{R}^n}u_1(x)\Phi(x)\,\mathrm{d} x+\frac{\varepsilon\beta}{2(\beta+1)}\left(1-\mathrm{e}^{-2t}\right)\int_{\mb{R}^n} u_{2}(x)\Phi(x)\,\mathrm{d}x\notag \\
& \geqslant \frac{\varepsilon}{2}\int_{\mb{R}^n}u_1(x)\Phi(x)\,\mathrm{d} x+\frac{\varepsilon\beta}{2(\beta+1)}\left(1-\mathrm{e}^{-1}\right)\int_{\mb{R}^n}u_2(x)\Phi(x)\,\mathrm{d}x\doteq C_2 \varepsilon
\end{align}
for any $t\geqslant1/2$. Here we remark that we may guarantee that $C_2>0$ because we assumed that at least one among the  nonnegative function $u_1$ or $u_2$ does not vanish identically.
%
%
%\color{blue}
%\noindent $\blacksquare$ 1. Although we may get lower bound for $F_1(t)\geqslant C_2\varepsilon$ for $t\geq1$, we should assume the nontriviality initial data $u_1$. Thus, I think the nontriviality initial data of $u_1$ or $u_2$ is relax. \textcolor{red}{I think you meant  $F_1(t)\geqslant C_2\varepsilon$ for $t\geqslant 0$. If it is so, I do agree with you of course, better to keep the assumptions on data as essential as possible}  \\
%\noindent $\blacksquare$ 2. In the next step, the first lower bound need to consider $t\geqslant1$. Therefore, it is not necessary for us to assume $t\geqslant0$. \textcolor{red}{Also this is a good remark, I agree} \\
%\noindent $\blacksquare$  3. In the next step, we need to use lower bound for $t\geqslant t_0$ to get iteration frame and $t_0\leqslant t/2\Rightarrow t\geqslant 2t_0$. Moreover, in the first lower bound, we assume $t\geqslant 1$. Thus, we consider $t_0=1/2$. \textcolor{red}{100 \% clear}
%
%\color{black}
By  H\"older's inequality and \eqref{L1 norm Psi B t+R}, we have
\begin{align*}
(1+\beta)H'(t)\geqslant C_1^{1-p} (t+R)^{-(n-1)(p-1)/2}(F_1(t))^p.
\end{align*}
Thus, integrating the above inequality over $[0,t]$ and using \eqref{Eq.05}, we obtain the iteration frame
\begin{align}\label{Iteration.Frame}
F_1(t)\geqslant C_3\int_0^t\mathrm{e}^{2(\tau-t)}\int_0^{\tau}(s+R)^{-(n-1)(p-1)/2}(F_1(s))^p\mathrm{d}s\,\mathrm{d}\tau,
\end{align}
where $C_3\doteq C_1^{1-p}/(1+\beta)$. We point out that in order to get \eqref{Iteration.Frame} we used the conditions $H(0)>0$ and  $F_1(0)\geqslant 0$.

The combination of \eqref{Initial.Lower.Bound} and \eqref{Iteration.Frame} shows %\textcolor{red}{$\blacktriangle$ here the steps are not really precise!}
%\begin{align*}
%\textcolor{red}{\blacktriangleright \mbox{old version} \ \ } F_1(t)&\geqslant C_2^pC_3 \, \varepsilon^p\int_{1/2}^t\mathrm{e}^{2(\tau-t)}\int_{1/2}^{\tau}(s+R)^{-(n-1)(p-1)/2}\mathrm{d}s\,\mathrm{d}\tau\\
%&\geqslant C_2^pC_3 \, \varepsilon^p(t+R)^{-(n-1)(p-1)/2}\int_{t/2}^t\mathrm{e}^{2(\tau-t)}\tau\mathrm{d}\tau\\
%&\geqslant 4^{-1} C_2^p C_3 \, \varepsilon^p(t+R)^{-(n-1)(p-1)/2}t\left(1-\mathrm{e}^{-t}\right).
%\end{align*}
\begin{align*}
%\textcolor{red}{\mbox{new version} \ \ }
F_1(t)&\geqslant C_2^pC_3 \, \varepsilon^p\int_{1/2}^t\mathrm{e}^{2(\tau-t)}\int_{1/2}^{\tau}(s+R)^{-(n-1)(p-1)/2}\mathrm{d}s\,\mathrm{d}\tau\\
&\geqslant C_2^pC_3 \, \varepsilon^p(t+R)^{-(n-1)(p-1)/2}\int_{t/2}^t\mathrm{e}^{2(\tau-t)}\left(\tau-1/2\right)\mathrm{d}\tau\\
&\geqslant 4^{-1} C_2^p C_3 \, \varepsilon^p(t+R)^{-(n-1)(p-1)/2}(t-1)\left(1-\mathrm{e}^{-t}\right).
\end{align*} for $t\geqslant 1$.
In particular, for $t\geqslant 1$ the factor containing the exponential function in the last line of the previous chain of inequalities can be estimate from below by a constant, namely,
\begin{align} \label{lower bound F 0 step}
F_1(t) \geqslant K_0 (t+R)^{-\alpha_0} \, (t-1)^{\gamma_0} \qquad \mbox{for any} \ t\geqslant 1,
\end{align} where the multiplicative constant is $K_0\doteq C_2^p C_3(1-\mathrm{e}^{-1}) \,  \varepsilon^p/4$ and the exponents are defined by $\alpha_0\doteq (n-1)(p-1)/2$ and $\gamma_0\doteq 1$.

\subsection{Iteration argument} \label{Subsection iteration argument subcrit}
The previous subsection is devoted to determine the iteration frame and a first lower bound for $F_1$. Our next goal is to derive a sequence of lower bounds for $F_1$ by using  \eqref{Iteration.Frame}. More precisely, we prove that
\begin{align}\label{sequence of lower bound F}
F_1(t) \geqslant K_j (t+R)^{-\alpha_j} (t- L_j)^{\gamma_j} \qquad \mbox{for any} \ t\geqslant L_j,
\end{align} where $\{K_j\}_{j\in \mathbb{N}}$, $\{\alpha_j\}_{j\in \mathbb{N}}$ and $\{\gamma_j\}_{j\in \mathbb{N}}$ are sequences of nonnegative real numbers that will be determined throughout this subsection and $\{L_j\}_{j\in\mathbb{N}}$ is the sequence of the partial products of the infinite product
\begin{align*}
\prod_{k=0}^\infty \ell_k \quad  \mbox{with} \ \ \ell_k\doteq 1+p^{-k} \ \ \mbox{for any} \ k\in \mathbb{N},
\end{align*} that is, $$L_j\doteq \prod_{k=0}^{j} \ell_k \quad \mbox{for any} \ j\in \mathbb{N}.$$ 
Clearly \eqref{lower bound F 0 step} implies \eqref{sequence of lower bound F} for $j=0$. We are going to show \eqref{sequence of lower bound F} via an inductive argument. Also, it remains to verify only the inductive step. Let us assume that \eqref{sequence of lower bound F} hols for $j\geqslant 0$. Then, in order to prove the inductive step, we shall prove \eqref{sequence of lower bound F} for $j+1$. After shrinking the domain of integration in \eqref{Iteration.Frame}, if we plug \eqref{sequence of lower bound F} in \eqref{Iteration.Frame}, we find
\begin{align*}
F_1(t)&\geqslant C_3K_j^p\int_{L_j}^t\mathrm{e}^{2(\tau-t)}\int_{L_j}^{\tau}(s+R)^{-(n-1)(p-1)/2-\alpha_jp}(s-L_j)^{\gamma_jp}\mathrm{d}s\,\mathrm{d}\tau\\
&\geqslant C_3K_j^p(t+R)^{-(n-1)(p-1)/2-\alpha_jp}\int_{L_j}^t\mathrm{e}^{2(\tau-t)}\int_{L_j}^{\tau}(s-L_j)^{\gamma_jp}\mathrm{d}s\,\mathrm{d}\tau\\
&\geqslant \frac{C_3K_j^p}{\gamma_jp+1}(t+R)^{-(n-1)(p-1)/2-\alpha_jp}\int_{t/\ell_{j+1}}^t\mathrm{e}^{2(\tau-t)}(\tau-L_j)^{\gamma_jp+1}\mathrm{d}\tau
\end{align*}
for any $t\geqslant L_{j+1} $. We point out that in the last step we could restrict the domain of integration with respect to $\tau$ from $[L_j,t]$ to $[t/\ell_{j+1},t]$ since $t\geqslant L_{j+1}$ and $\ell_{j+1}>1$ imply the inequality $L_{j} \leqslant t/\ell_{j+1}<t$. Also,
\begin{align*}
F_1(t)\geqslant\frac{C_3K_j^p}{2(\gamma_jp+1)\ell_{j+1}^{\gamma_jp+1}}(t+R)^{-(n-1)(p-1)/2-\alpha_jp}(t-L_{j+1})^{\gamma_jp+1}\left(1-\mathrm{e}^{2t(1/\ell_{j+1}-1)}\right)
\end{align*}
for any $t\geqslant L_{j+1} $. We observe that for $t\geqslant L_{j+1} \geqslant \ell_{j+1} $ it is possible to estimate 
\begin{align}\label{Ineq.1-e}
1-\mathrm{e}^{2t(1/\ell_{j+1}-1)} & \geqslant 1-\mathrm{e}^{-2(\ell_{j+1}-1)}\geqslant  2(\ell_{j+1}-1)(2-\ell_{j+1})\notag\\
& \geqslant 2(p^{j+1}-1)p^{-2(j+1)}\geqslant 2(p-1)p^{-2(j+1)}.
\end{align}
Thus, for any $t\geqslant L_{j+1}$ we have proved
\begin{align*}
F_1(t)  & \geqslant  \frac{(p-1)p^{-2(j+1)}C_3K_j^p}{(\gamma_jp+1)\ell_{j+1}^{\gamma_jp+1}} (t+R)^{-(n-1)(p-1)/2-\alpha_jp}(t-L_{j+1})^{\gamma_jp+1},
\end{align*}
which is exactly \eqref{sequence of lower bound F} for $j+1$, provided that
\begin{align*}
K_{j+1} & \doteq   \frac{(p-1)p^{-2(j+1)}C_3K_j^p}{(\gamma_jp+1)\ell_{j+1}^{\gamma_jp+1}}, \ \
\alpha_{j+1}  \doteq \frac{1}{2}(n-1)(p-1)+\alpha_jp  , \ \
\gamma_{j+1}  \doteq  \gamma_jp+1 .
\end{align*}
By using recursively the previous relations for $\alpha_j$ and $\gamma_j$ it is easy to get
\begin{align*}
\alpha_j&=p^j\left(\alpha_0+\tfrac{n-1}{2}\right)-\tfrac{n-1}{2},\\
\gamma_j&=p^j\left(\gamma_0+\tfrac{1}{p-1}\right)-\tfrac{1}{p-1}.
\end{align*}  Besides, the inequality $\gamma_{j-1}p+1= \gamma_j \leqslant p^j\left(\gamma_0+\tfrac{1}{p-1}\right)$ implies immediately
\begin{align*}
K_j&\geqslant(p-1)C_3\left(\gamma_0+\tfrac{1}{p-1}\right)^{-1}K_{j-1}^pp^{-3j}\ell_j^{-\gamma_j}.
\end{align*}
 Due to the choice of $\ell_j$, it holds
\begin{align*}
\lim_{j\to \infty} \ell_j^{\gamma_j}= \lim_{j\to \infty} \exp\left(\left(\gamma_0+\tfrac{1}{p-1}\right)p^{j}\log \left(1+p^{-j}\right)  \right) = \mathrm{e}^{\gamma_0+1/(p-1)}.
\end{align*} Therefore, there exists a suitable constant $M=M(n,p)>0$ such that $\ell_j^{-\gamma_j}\geqslant M$ for any $j\in \mathbb{N}$. So, combining this inequality with the previous estimate from below of $K_j$, we have
\begin{align*}
K_j \geqslant \underbrace{ (p-1)MC_3\left(\gamma_0+\tfrac{1}{p-1}\right)^{-1}}_{\doteq D}K_{j-1}^p p^{-3j} \quad \mbox{for any} \ j\in\mathbb{N}.
\end{align*}
If we apply the logarithmic function to both sides of the inequality $K_j\geqslant D K_{j-1}^p p^{-3j}$ and we use iteratively the resulting inequality, we obtain
\begin{align*}
\log K_j&\geqslant p^{j}\log K_0 -3 \left(\sum_{k=0}^{j-1} (j-k)p^k\right) \log p+ \left(\sum_{k=0}^{j-1} p^k\right)\log D\\
&\geqslant p^j \left(\log K_0-\frac{3p \log p}{(p-1)^2} +\frac{\log D}{p-1} \right)+\frac{3j \log p}{p-1}+\frac{3p \log p}{(p-1)^2}-\frac{\log D}{p-1}
\end{align*}
for any $j\in\mathbb{N}$, where in the second step we use the identity $$\sum_{k=0}^{j-1} (j-k)p^k=\frac{1}{p-1}\left(\frac{p^{j+1}-p}{p-1}-1\right).$$ Let $j_0=j_0(n,p)\in\mathbb{N}$ be the smallest nonnegative integer such that
\begin{align*}
j_0\geqslant \frac{\log D}{3\log p}-\frac{p}{p-1}.
\end{align*}
 Then, for any $j\geqslant j_0$ it results
 \begin{align*}\label{lower bound log Kj}
 \log K_j\geqslant p^j\log\left(D^{1/(p-1)}p^{-3p/(p-1)^2}K_0\right)=p^j\log(E\varepsilon^p)
 \end{align*}
 for a suitable positive constant $E=E(n,p)$.
 
  Let us denote $$L\doteq \lim_{j\to\infty} L_j = \prod_{j=0}^\infty \ell_{j}\in \mathbb{R}.$$  Thanks to $\ell_j>1$, it holds $L_j \uparrow L$ as $j\to \infty$. In particular, \eqref{sequence of lower bound F} holds for any $j\in\mathbb{N}$ and any $t\geqslant L$.
 
  Combining the above results and using the explicit representation for $\alpha_j$ and $\gamma_j$, we get
 \begin{align*}
 F_1(t)&\geqslant \exp\left(p^j\log(E\varepsilon^p)\right)(t+R)^{-\alpha_j}(t-L)^{\gamma_j}\\
 &\geqslant \exp\left(p^j\left(\log(E\varepsilon^p)-\left(\alpha_0+\tfrac{n-1}{2}\right)\log(t+R)+\left(\gamma_0+\tfrac{1}{p-1}\right)\log(t-L)\right)\right)\\
 &\quad\times(t+R)^{(n-1)/2}(t-L)^{-1/(p-1)}
 \end{align*}
 for any $j\geqslant j_0$ and any $t\geqslant L$.
 
 Then, since for $t\geqslant \max\{R,2L \}$ we may estimate $R+t\leqslant 2t$ and $t-L \geqslant t/2$, we find
\begin{equation}\label{final lower bound F_0}
F_1(t)  \geqslant  \exp \left(p^j \log \left(E_1\varepsilon^p t^{\gamma_0+\frac{1}{p-1}-(\alpha_0+\frac{n-1}{2})}\right)\right) (t+R)^{n}(t-L)^{-1/(p-1)}
\end{equation}
for any $j\geqslant j_0$, where $E_1\doteq2^{-\left(\alpha_0+(n-1)/2+\gamma_0+1/(p-1)\right)}E$.
 We rewrite the exponent for $t$ in the last inequality as follows:
 \begin{align*}
 \gamma_0+\tfrac{1}{p-1}-\left(\alpha_0+\tfrac{n-1}{2}\right)=\tfrac{p}{2(p-1)}((n+1)-(n-1)p)&=\tfrac{p((n+1)-(n-1)p)}{2(p-1)}.
 \end{align*}
 We notice that for $1<p<p_{\mathrm{Gla}}(n)$ (respectively, for $1<p$ when $n=1$), this exponent for $t$ is positive. Let us fix $\varepsilon_0=\varepsilon_0(u_0,u_1,u_2,n,p,R,\beta)>0$ such that
 \begin{align*}
 \varepsilon_0^{-\frac{2(p-1)}{(n+1)-(n-1)p}}\geqslant E_1^{\frac{2(p-1)}{p((n+1)-(n-1)p)}} \max\{R,2L\}.
 \end{align*} Also, for any $\varepsilon\in(0,\varepsilon_0]$ and any $t>E_2  \varepsilon^{-\frac{2(p-1)}{(n+1)-(n-1)p}}$, where $E_2\doteq E_1^{-\frac{2(p-1)}{p((n+1)-(n-1)p)}}$, we have
 \begin{align*}
 t\geqslant \max\{R,2L\} \quad \mbox{and} \quad \log \bigg(E_1\varepsilon^p t^{ \tfrac{p((n+1)-(n-1)p)}{2(p-1)}}\bigg)>0.
 \end{align*} Consequently, for any $\varepsilon\in(0,\varepsilon_0]$ and any $t>E_2  \varepsilon^{-\frac{2(p-1)}{(n+1)-(n-1)p}}$ letting $j\to \infty$ in \eqref{final lower bound F_0} we find that the lower bound for $F_1$ blows up. So, for any $\varepsilon\in(0,\varepsilon_0]$ the functional $F_1$ has to blow up in finite time too and, furthermore, the lifespan of the local solution $u$ can be estimated from above in the following way: $$T(\varepsilon) \leqslant C \varepsilon^{-(\frac{1}{p-1} -\frac{n-1}{2})^{-1}}.$$ We completed the proof of Theorem \ref{Thm blow-up} in the case $1<p<p_{\mathrm{Gla}}(n)$.  In the next section we will investigate the blow -- up dynamic in the case $p= p_{\mathrm{Gla}}(n)$.

 \section{Blow -- up result in the critical case}\label{Section.Blow-up.Critical}

\subsection{Iteration frame}
From the last section, we know that the first lower bound for functional $F_1$ is given by
\begin{align*}
F_1(t)\geqslant C_2\varepsilon
\end{align*}
for any $t\geqslant1/2$, with a positive constant $C_2$.

In this section, we consider the case $p=p_{\text{Gla}}(n)=(n+1)/(n-1)$ when $n\geqslant 2$. In this special case, the iteration frame  \eqref{Iteration.Frame} takes the form
\begin{align}\label{Iteration.Frame2}
F_1(t)&\geqslant C_3\int_{0}^t\mathrm{e}^{2(\tau-t)}\int_{0}^{\tau} (s+R)^{-1}(F_1(s))^p\mathrm{d}s\,\mathrm{d}\tau\notag\\
&\geqslant C_4\int_{1}^t\mathrm{e}^{2(\tau-t)}\int_{1}^{\tau}\frac{(F_1(s))^p}{s}\mathrm{d}s\,\mathrm{d}\tau
\end{align}
for some suitable positive constant $C_4$ and for any $t\geqslant1$.

\subsection{Iteration argument} \label{Subsection iteration argument crit}
 Analogously to what we did in Subsection \ref{Subsection iteration argument subcrit}  we derive now a sequence of lower bounds for $F_1$ by using the iteration frame \eqref{Iteration.Frame2}. More specifically, we want to show that
\begin{align}\label{sequence of lower bound F2}
F_1(t) \geqslant Q_j \left(\log(t/L_j)\right)^{\sigma_j} \qquad \mbox{for any} \ t\geqslant L_j ,
\end{align} where $\{Q_j\}_{j\in \mathbb{N}}$ and $\{\sigma_j\}_{j\in \mathbb{N}}$  are sequences of nonnegative real numbers to be determined and $\{L_j\}_{j\in\mathbb{N}}$ is defined as in Section \ref{Section.Blow-up.Subcritical}. When $j=0$, we have $Q_0\doteq C_2\varepsilon$ and $\sigma_0\doteq 0$ according to \eqref{Initial.Lower.Bound}.

As in the subcritical case, we are going to prove \eqref{sequence of lower bound F2} by using an inductive argument.  We assume the validity of \eqref{sequence of lower bound F2} for $j\geqslant 0$ and we have to prove it for $j+1$, prescribing the values of $Q_{j+1}$ and of $\sigma_{j+1}$. Shrinking the domain of integration in \eqref{Iteration.Frame2} and plugging \eqref{sequence of lower bound F2} in \eqref{Iteration.Frame2}, we obtain
\begin{align*}
F_1(t)&\geqslant C_4Q_j^p\int_{L_j}^t\mathrm{e}^{2(\tau-t)}\int_{L_j}^{\tau}\frac{(\log(s/L_j))^{\sigma_jp}}{s}\mathrm{d}s\,\mathrm{d}\tau\\
&\geqslant \frac{C_4Q_j^p}{\sigma_jp+1}\int_{L_j}^t\mathrm{e}^{2(\tau-t)}(\log(\tau/L_j))^{\sigma_jp+1}\mathrm{d}\tau
\end{align*}
for any $t\geqslant L_{j+1}$. Since for $t\geqslant L_{j+1}$ it holds $L_j\leqslant t/\ell_{j+1}$, a restriction of the domain of integration in the last inequality yields
\begin{align*}
F_1(t)& \geqslant\frac{C_4Q_j^p}{\sigma_jp+1}\int_{t/\ell_{j+1}}^t\mathrm{e}^{2(\tau-t)}(\log(\tau/L_j))^{\sigma_jp+1}\mathrm{d}\tau\\
&\geqslant\frac{C_4Q_j^p}{\sigma_jp+1}(\log(t/L_{j+1}))^{\sigma_jp+1}\int_{t/\ell_{j+1}}^t\mathrm{e}^{2(\tau-t)}\mathrm{d}\tau\\
&\geqslant\frac{C_4Q_j^p}{2(\sigma_jp+1)}(\log(t/L_{j+1}))^{\sigma_jp+1}\left(1-e^{2t(1/\ell_{j+1}-1)}\right)\\
&\geqslant\frac{C_4Q_j^p(p-1)p^{-2(j+1)}}{\sigma_jp+1}(\log(t/L_{j+1}))^{\sigma_jp+1},
\end{align*}
where we used once again \eqref{Ineq.1-e} in the last inequality. So, we proved \eqref{sequence of lower bound F2} for $j+1$, provided that
\begin{align*}
Q_{j+1} & \doteq   \frac{C_4Q_j^p(p-1)p^{-2(j+1)}}{\sigma_jp+1}, \ \
\sigma_{j+1}  \doteq \sigma_jp+1 .
\end{align*}
Repeating the same procedure seen in Section \ref{Section.Blow-up.Subcritical}, we get easily
\begin{align*}
\sigma_j&=\tfrac{p^j-1}{p-1},\\
Q_j&\geqslant C_4(p-1)^2p^{-3j}Q_{j-1}^p\doteq \widetilde{D}p^{-3j}Q_{j-1}^p.
\end{align*}
 Hence, applying again the monotonicity of the logarithmic function, in this case to the inequality $Q_j \geqslant \widetilde{D}p^{-3j}Q_{j-1}^p$, we derive 
\begin{align*}
\log Q_j&\geqslant p^j \left(\log Q_0-\frac{3p \log p}{(p-1)^2} +\frac{\log \widetilde{D}}{p-1} \right)+\frac{3j \log p}{p-1}+\frac{3p \log p}{(p-1)^2}-\frac{\log \widetilde{D}}{p-1}
\end{align*}
for any $j\in\mathbb{N}$. Let $j_1=j_1(n,p)\in\mathbb{N}$ be the smallest nonnegative integer such that
\begin{align*}
j_1\geqslant \frac{\log \widetilde{D}}{3\log p}-\frac{p}{p-1}.
\end{align*}
Then, for any $j\geqslant j_1$ it results
\begin{align*}\label{lower bound log Kj2}
\log Q_j\geqslant p^j\log\left(\widetilde{D}^{1/(p-1)}p^{-3p/(p-1)^2}Q_0\right)=p^j\log(\widetilde{E}\varepsilon)
\end{align*}
for a suitable positive constant $\widetilde{E}=\widetilde{E}(n,p)$.
 Let us recall that $L$ denotes the monotonic limit of the sequence $\{L_j\}_{j\in\mathbb{N}}$. Therefore, we have that \eqref{sequence of lower bound F2} holds for any $j\in\mathbb{N}$ and any $t\geqslant L$.

Thus, applying the explicit representation for $\sigma_j$, we arrive at
\begin{equation}\label{Final Lower bound 2}
\begin{aligned}
F_1(t)&\geqslant \exp\left(p^j\log(\widetilde{E}\varepsilon)\right)(\log(t/L))^{\sigma_j}\\
&=\exp\left(p^j\log\left(\widetilde{E}\varepsilon(\log(t/L))^{1/(p-1)}\right)\right)(\log(t/L))^{-1/(p-1)},
\end{aligned}
\end{equation} for any $ j\geqslant j_1$ and any $t\geqslant L$.
 In this case, we fix $\varepsilon_0=\varepsilon_0(u_0,u_1,u_2,n,p,R,\beta)>0$ in such a way that
\begin{align*}
\exp\left(\widetilde{E}^{-p+1}\varepsilon_0^{-(p-1)}\right)\geqslant1.
\end{align*}
 Consequently, for any $\varepsilon\in(0,\varepsilon_0]$ and any $t>L\exp(\widetilde{E}^{-p+1}\varepsilon^{-(p-1)}) $, we get
\begin{align*} t\geqslant L \quad \mbox{and} \quad \log \bigg(\widetilde{E}\varepsilon(\log(t/L))^{1/(p-1)}\bigg)>0.
\end{align*} 
Hence, for any $\varepsilon\in(0,\varepsilon_0]$ and any $t>L\exp(\widetilde{E}^{-p+1}\varepsilon^{-(p-1)})$ by letting $j\to \infty$ in \eqref{Final Lower bound 2} we see that the lower bound for $F_1$ blows up. Thus, for any $\varepsilon\in(0,\varepsilon_0]$ the functional $F_1$ has to blow up in finite time as well and, besides, the lifespan of the local solution $u$ can be estimated from above in the following way $$T(\varepsilon) \leqslant \exp\left(C\varepsilon^{-(p-1)}\right),$$ for a suitable constant $C$ which is independent of $\varepsilon$.   This completes the proof of Theorem \ref{Thm blow-up} in the case $p=p_{\mathrm{Gla}}(n)$.

\section{ Final remarks }

In Theorem \ref{Thm blow-up}, we proved a blow -- up result for $1<p\leqslant p_{\mathrm{Gla}}(n)$ under suitable sign and support assumptions for the Cauchy data. Furthermore, as byproduct of the iteration arguments we obtained upper bound estimates for the lifespan as well. In particular, we find the same range for $p$ in the blow -- up result as for the semilinear Cauchy problem \eqref{Semi.Wave.Eq}, which is known to be sharp in the case of this last wave model.

\subsection*{Acknowledgments} The Ph.D. study of the first author is supported by S\"achsiches Landesgraduiertenstipendium. The second author is supported by the University of Pisa, Project PRA 2018 49.

%\medskip

% ------------------------------------------------------------------------
\end{document}